\documentclass{amsart}%
\usepackage{amsmath}
\usepackage{graphicx}%
\usepackage{amsfonts}%
\usepackage{amssymb}

\theoremstyle{plain}

\newtheorem{corollary}{Corollary}

\newtheorem{lemma}{Lemma}

\newtheorem{proposition}{Proposition}

\numberwithin{equation}{section}

\begin{document}
\title[Semi-infinite VHS and integrable hierarchies.]{Semi-infinite variations of Hodge structures and integrable hierarchies of KdV-type}
\author{Serguei Barannikov}
\address{Ecole Normale Superieure, 45 rue d'Ulm, Paris 75230}
\email{serguei.barannikov@ens.fr}
\urladdr{http://www.dma.ens.fr/\symbol{126}barannik}
\date{April 29, revised August 21}

\begin{abstract}
We introduce integrable KdV-type hierarchy associated naturally with arbitrary
semi-simple Frobenius manifold. We present hierarchy in a Lax form and show
that it admits bihamiltonian description. The hierarchy allows to extend
corresponding semi-infinite variation of Hodge structures by including
variation along higher times $\{x^{\alpha,r}\}$ satisfying $\frac{\partial
}{\partial x^{\alpha,r}}\mathcal{L}(x)\subseteq\hbar^{-r}\mathcal{L}(x)$

\end{abstract}
\maketitle

\section{Introduction}

Geometry of families of semi-infinite subspaces plays central role in elegant
approach to integrable hierarchies in \cite{SW}. Similar geometry can be used
to describe Frobenius manifolds (see \cite{B1}, \cite{B2}). Namely there is
canonical family of Frobenius manifold structures on parameter space of any
abstract semi-infinite variation of Hodge structures of Calabi-Yau type and
conversely one can associate such semi-infinite variation of Hodge structures
($\frac{\infty}{2}-$VHS for short) with arbitrary Frobenius manifold (we
recall this in section 2 below).   In this note we pursue analysis of
relations between $\frac{\infty}{2}-$VHS and integrable hierarchies and
introduce integrable KdV type hierarchy which is associated with $\frac
{\infty}{2}-$VHS corresponding to an arbitrary \emph{semisimple} Frobenius
manifold. This note arose from an attempt to understand relation between
problem of extension of abstract $\frac{\infty}{2}-$VHS of Calabi-Yau type and
dressing method from theory of integrable hierarchies (see \cite{DS} \S 1 and
references therein).

\subsection{Some notations:}

$Mat(n,\mathbb{C})$ denotes algebra of $(n\times n)$ matrices, $Diag\subset
Mat(n,\mathbb{C})$ denotes subalgebra of matrices whose entries are zero apart
from diagonal, for $A\in Mat(n,\mathbb{C})$ we denote via $A^{\top}$
transposed matrix, for a vector space $V$ we denote via $V((\hbar))$ the space
of Laurent series with values in $V$, for $v\in V((\hbar))$, $v=\sum
_{l=-m}^{+\infty}v_{l}\hbar^{l}$ we set $v_{\geq k}=\sum_{l=k}^{+\infty}%
v_{l}\hbar^{l},v_{<k}=\sum_{l=-m}^{l=k-1}v_{l}\hbar^{l}$.

\section{\bigskip Semi-infinite VHS and Frobenius manifolds.}

In this section I recall construction of family of Frobenius manifolds
associated with abstract $\frac{\infty}{2}-$VHS of Calabi-Yau type described
in \cite{B1} and give also inverse construction. In particular it provides a
nice mathematical explanation to universality of WDVV-equations. In this
section we work in analytic category and pure even case leaving appropriate
adjustments for other categories (formal, $\mathbb{Z}/2\mathbb{Z}-$graded
etc.) to an interested reader.

Let $Gr_{\frac{\infty}{2}}^{(n)}$denotes affine grassmanian (see \cite{PS},
\S\ 8). Recall that it is defined as quotient $Gr_{\frac{\infty}{2}}%
^{(n)}:=LGL(n,\mathbb{C})/L^{+}GL(n,\mathbb{C})$ where $LGL(n,\mathbb{C})$
denotes group of maps from circle $S^{1}=\{\hbar\in\mathbb{C}:\,|\hbar|=1\}$
to $GL(n,\mathbb{C})$ which are analytic in some neighborhood of $S^{1}$ and
$L^{+}GL(n,\mathbb{C})$ denotes subgroup of elements which are boundary values
of analytic maps from disk $\{\hbar\in\mathbb{C}:\,|\hbar|\leq1\}$ to
$GL(n,\mathbb{C})$.

Let $H^{(n)}$ denotes Hilbert space of all square-integrable functions on
circle $S^{1}=\{\hbar\in\mathbb{C}:\,|\hbar|=1\}$ with values in
$\mathbb{C}^{n}$, $H^{(n)}=L^{2}(S^{1},\mathbb{C}^{n})$. Let us denote also
via $H_{+}^{(n)}$ (resp. $H_{-}^{(n)}$) closed subspace of $H^{(n)}$ generated
by elements of form $v\cdot\hbar^{k}$, $v\in\mathbb{C}^{n}$, $k\geq0$ (resp.
$k<0$), so that $H^{(n)}=H_{+}^{(n)}\oplus H_{-}^{(n)}$, and by $pr_{+}$,
(resp. $pr_{-}$) the orthogonal projection $H^{(n)}\rightarrow$ $H_{+}^{(n)}$
(resp. $H^{(n)}\rightarrow$ $H_{-}^{(n)}$) along $H_{-}^{(n)}$ (resp.
$H_{+}^{(n)}$).

\begin{lemma}
(\cite{PS},\S 8) Grassmanian $Gr_{\frac{\infty}{2}}^{(n)}$ can be defined
alternatively as set of all closed subspaces $\mathcal{L}\subset$ $H^{(n)}$
having everywhere dense subset consisting of analytic functions and such that:

\begin{itemize}
\item $pr_{+}|_{\mathcal{L}}$ is a Fredholm operator (recall that operator $T$
is Fredholm iff $\dim$ ker $T$, $\dim$ coker $T<\infty$), i.e. $\mathcal{L}$
is in a sense ''comparable'' with $H_{+}^{(n)}$,

\item $\hbar\mathcal{L}\subset\mathcal{L}$.
\end{itemize}
\end{lemma}

\begin{proof}
To a class $[\varphi]\in LGL(n,\mathbb{C})/L^{+}GL(n,\mathbb{C})$ one can
associate subspace $\varphi\cdot H_{+}^{(n)}\subset H^{(n)}$. The Fredholm
property follows from (\cite{PS}, proposition 6.3.1). Conversely, the
factorspace $\mathcal{L}/\hbar\mathcal{L}$ \ for a subspace $\mathcal{L}$
satisfying the above properties is an $n-$dimensional vector space, since
inclusion $\hbar\mathcal{L}\subset\mathcal{L}$ is a Fredholm operator of index
equal to the index of inclusion $\hbar H_{+}^{(n)}\subset H_{+}^{(n)},$ i.e.
$n$; therefore if $\varphi_{i}\in\mathcal{L}$, $i=1,\ldots,n$, are such that
$\{\varphi_{i}$ mod $\hbar\mathcal{L}\}_{i=1\ldots n}$ is a basis for
$\mathcal{L}/\hbar\mathcal{L}$, then matrix with columns $\varphi_{i}$ defines
the corresponding element from quotient $LGL(n,\mathbb{C})/L^{+}%
GL(n,\mathbb{C})$.
\end{proof}

Let $\mathcal{L}(x)\in Gr_{\frac{\infty}{2}}^{(n)}$, $x\in\mathcal{U}$, be a
family of subspaces from $Gr_{\frac{\infty}{2}}^{(n)}$ parametrized by
$\mathcal{U}$.

Our first assumption on family $\mathcal{L}(x)$ is

\emph{1) }$\frac{\infty}{2}-$\emph{Griffiths transversality: }$\frac{\partial
}{\partial x}\mathcal{L}(x)\subseteq\hbar^{-1}\mathcal{L}(x)$

For family of subspaces having such property one has ''symbol of
$\frac{\partial}{\partial x}$'' map:
\[
Symbol(\frac{\partial}{\partial x}):T_{x}\mathcal{U}\otimes\mathcal{L}%
/\hbar\mathcal{L}\rightarrow\hbar^{-1}\mathcal{L}/\mathcal{L}%
\]
where $T_{x}\mathcal{U}$ denotes tangent space to $\mathcal{U}$ at a point $x
$. Our next assumption on family $\mathcal{L}(x)$ is

\emph{2)''Calabi-Yau type''} : there exists one-dimensional subspace
$\{\lambda\lbrack\Omega]\mathbb{\}}_{\lambda\in\mathbb{C}}\subset
\mathcal{L}/\hbar\mathcal{L}(x)$ for any $x\in\mathcal{U}$, such that map
$Symbol(\frac{\partial}{\partial x})|_{\{\lambda\lbrack\Omega]\}}%
:T_{x}\mathcal{U}\otimes\{\lambda\lbrack\Omega]\}\rightarrow\hbar
^{-1}\mathcal{L}/\mathcal{L} $ is an isomorphism.

In particular, $\dim_{\mathbb{C}}\mathcal{U}=\dim_{\mathbb{C}}\mathcal{L}%
/\hbar\mathcal{L}$ and also the map $T_{x}\mathcal{U}\rightarrow
Hom(\mathcal{L}/\hbar\mathcal{L},\hbar^{-1}\mathcal{L}/\mathcal{L)}$ induced
by $Symbol(\frac{\partial}{\partial x})$ is an embedding. Next we need to
incorporate a semi-infinite analog of Poincare pairing into our assumptions.
Let $G:(H^{(n)})^{\otimes2}\rightarrow H^{(1)}$ be a linear nondegenerate
pairing, which is symmetric in the following sense: $G(a,b)(\hbar
)=(-1)^{n}G(b,a)(-\hbar)$ and has the following property of linearity with
respect to multiplication by $\hbar$ : $G(\hbar a,b)=G(a,-\hbar b)=\hbar
G(a,b)$. Such pairing is uniquely defined by restriction $G|_{(\hbar^{0}%
\cdot\mathbb{C}^{n})^{\otimes2}}=\sum_{-\infty}^{+\infty}\hbar^{i}g^{(i)}$
where $g^{(i)}:(\mathbb{C}^{n})^{\otimes2}\rightarrow\mathbb{C}$ is
$(-1)^{n+i}-$symmetric. Our last assumption on family $\mathcal{L}(x)$ is

\emph{3) Isotropy with respect to pairing:} $G|_{\mathcal{L}^{\otimes2}}%
\in\hbar^{n}H_{+}^{(1)}$

Examples of families of $\frac{\infty}{2}$-subspaces satisfying these three
conditions arise naturally in context of noncommutative algebraic geometry
(see \cite{B1} where also relation with standard variations of Hodge
structures is explained). In a sense one can say that such family of
semi-infinite subspaces indicates presence of a non-commutative complex
Calabi-Yau manifold.

Let $Gr_{-\frac{\infty}{2}}^{(n)}$ denotes ''opposite'' grassmanian consisting
of closed subspaces of $H^{(n)}$ such that for $S\in$ $Gr_{-\frac{\infty}{2}%
}^{(n)}$ restriction on $S$ of projection $pr_{-}|_{S}$ to $H_{-}^{(n)}$ is a
Fredholm operator, that analytic functions are dense in $S$ and that
$\hbar^{-1}S\subset S$. Let $S\in Gr_{-\frac{\infty}{2}}^{(n)}$ be a subspace
which is transversal to $\mathcal{L}(x)$ for all $x\in\mathcal{U}$:
$S\oplus\mathcal{L}(x)=H^{(n)}$ and which satisfies in addition the following
isotropy condition $G|_{S^{\otimes2}}\in\hbar^{n-1}H_{-}^{(1)}$.
Transversality implies that intersection $\mathcal{L}(x)\cap\hbar S$ is an
$n-$dimensional vector space. There are natural isomorphisms: $i_{S}%
:\mathcal{L}(x)\cap\hbar S$ $\simeq\hbar S/S$ and $i_{\mathcal{L}}%
:\mathcal{L}(x)\cap\hbar S$ $\simeq\mathcal{L}/\hbar\mathcal{L}$. Let
$\omega\in\hbar S/S$. Constraint imposed on element $[\Omega]$ from the second
condition is of an open type. Setting if necessary $\mathcal{U}$ to be its
open subset we can assume that there exists $\omega$ such that restriction to
$i_{\mathcal{L}(x)}(i_{S}^{-1}\omega)$ of $Symbol(\frac{\partial}{\partial
x})$ is an isomorphism for all $x\in\mathcal{U}$. Let $\psi(x)\in
\mathcal{L}(x)\cap\hbar S$ be the unique element such that $i_{S}%
(\psi(x))=\omega$. It is given by intersection of $\mathcal{L}(x)$ with
constant affine space $S+\omega$. Let us consider map $\mathcal{U}\rightarrow
S/\hbar^{-1}S$ which sends $x$ to class of $[\psi(x)-\omega]$. It follows that
differential of this map is an isomorphism. Let us choose a basis
$\{\Delta_{\alpha}\}_{\alpha\in\lbrack1,\ldots,n]}$ in $S/\hbar^{-1}S$ so that
$\Delta_{1}=\hbar^{-1}\omega$, and denote via $\{x^{\alpha}\}$ the
corresponding coordinates on $\mathcal{U}$ induced from linear coordinates on
$S/\hbar^{-1}S$ via the map $[\psi(x)-\omega]$. Then $\{\frac{\partial\psi
(x)}{\partial x^{\alpha}}\}_{\alpha\in\lbrack1,\ldots,n]}$ is a basis in
$\hbar^{-1}\mathcal{L}(x)\cap S$ and $\mathcal{L}(x)=[\frac{\partial\psi
(x)}{\partial x}]\hbar H_{+}^{(n)}$.

\begin{proposition}
Element $\psi(x)$ satisfy
\begin{align}
\frac{\partial^{2}\psi(x)}{\partial x^{\alpha}\partial x^{\beta}}  &
=\hbar^{-1}\sum_{\gamma}C_{\alpha\beta}^{\gamma}(x)\frac{\partial\psi
(x)}{\partial x^{\gamma}}\label{eqprfrob}\\
G(\frac{\partial\psi(x)}{\partial x^{\alpha}},\frac{\partial\psi(x)}{\partial
x^{\beta}})  &  =\hbar^{n-2}\eta_{\alpha\beta},\,\,\,\,\eta_{\alpha\beta
}=const\nonumber\\
\frac{\partial\psi(x)}{\partial x^{1}}  &  =\hbar^{-1}\omega\nonumber
\end{align}
\end{proposition}

\begin{proof}
See proofs of propositions 6.5 from \cite{B1}, 4.1, 4.4, 4.8 from \cite{B2}.
\end{proof}

Now same arguments as in  \cite{B2}, \S 4 give the following corollary (we
refer reader to \cite{M} for a definition of Frobenius manifold)

\begin{corollary}
\label{corfrob}Tensors $C_{\alpha\beta}^{\gamma}(x^{\alpha}),\,\eta
_{\alpha\beta},\frac{\partial}{\partial x^{1}}$ define Frobenius manifold
structure on $\mathcal{U}$. \thinspace
\end{corollary}

Remark that conformal (i.e. equipped with Euler vector field) Frobenius
manifold corresponds to the data as above equipped with constant first order
differential operator $\mathcal{D}_{\partial/\partial\hbar}$ acting on
elements of $H^{(n)}$ such that $\mathcal{D}_{\partial/\partial\hbar
}(\mathcal{L})\subseteq\hbar^{-2}\mathcal{L}$, $\mathcal{D}_{\partial
/\partial\hbar}(S)\subseteq\hbar^{-1}S$ and $\mathcal{D}_{\partial
/\partial\hbar}G=0$.

Conversely let we are given Frobenius manifold structure $(C_{\alpha\beta
}^{\gamma}(x^{\alpha})_{x\in\mathcal{U}},$ $\eta_{\alpha\beta},$
$e=\frac{\partial}{\partial x^{1}})$ on $\mathcal{U}$. Then operators
$\frac{\partial}{\partial x^{\alpha}}-\hbar^{-1}C_{\alpha}(x)$, $\alpha
\in\lbrack1,\ldots,n]$, commute with each other and therefore there exists
fundamental solution $\phi_{\alpha}^{i}$, $i\in\lbrack1,\ldots,n]\,$, to
\begin{equation}
\frac{\partial}{\partial x^{\alpha}}\phi_{\beta}(x,\hbar)=\hbar^{-1}%
\sum_{\gamma}C_{\alpha\beta}^{\gamma}(x)\phi_{\beta}(x,\hbar),\,\,\,\alpha
,\beta\in\lbrack1,\ldots,n]\,\,\,\label{dpsi}%
\end{equation}
so that $\det\phi_{\alpha}^{i}(x,\hbar)\neq0$. Such solution is unique up to
multiplication by an element from $LGL(n,\mathbb{C})$. We have $\hbar^{-1}%
\phi_{\alpha}^{i}=\frac{\partial}{\partial x^{\alpha}}\phi_{1}^{i}$. We define
corresponding family of semi-infinite subspaces as $\mathcal{L}(x):=\phi
(x)H_{+}^{(n)}$, i.e. it is the subspace generated by vectors $\hbar^{k}%
\phi_{\alpha}$, $\phi_{\alpha}=\sum_{i}\phi_{\alpha}^{i}\rho_{i}$ where
$\rho_{i}$ is the standard basis in $\mathbb{C}^{n}$ and $\alpha\in
\lbrack1,\ldots,n]$, $k\geq0$. It follows from (\ref{dpsi}) that it satisfies
$\frac{\infty}{2}-$Griffiths transversality condition. If $Symbol(\frac
{\partial}{\partial x})$ is written using basis $\frac{\partial}{\partial
x^{\alpha}},[\phi_{\alpha}],[\hbar^{-1}\phi_{\alpha}]$ in $T_{x}\mathcal{U}$,
$\mathcal{L}/\hbar\mathcal{L}$, $\hbar^{-1}\mathcal{L}/\mathcal{L}$
respectively then $Symbol(\frac{\partial}{\partial x})$ coincides with tensor
$C_{\alpha\beta}^{\gamma}(x)$ of multiplication. Therefore the second
condition is satisfied as one can take for $[\Omega]$ any element
corresponding to an invertible element in the algebra defined by
$C_{\alpha\beta}^{\gamma}$. Next we define pairing $G(a,b):=\hbar^{n-2}%
\eta(\hbar\phi^{-1}a,\hbar\phi^{-1}b)$ where we set $\eta(\hbar^{k}v,\hbar
^{l}u):=\hbar^{k}(-\hbar)^{l}\sum_{\alpha,\beta}$ $\eta_{\alpha\beta}%
v^{\alpha}u^{\beta}$ for $v=(v^{1},\ldots,v^{n})$, $u=(u^{1},\ldots,u^{n})$.
It follows from definition that $G|_{\mathcal{L}^{\otimes2}}\in\hbar^{n}%
H_{+}^{(1)}$. Notice that because of compatibility of $\eta$ with
multiplication defined by $C_{\alpha\beta}^{\gamma}(x)$ the value of $G(a,b)$
does not depend on $x$. It follows that for any choice of opposite isotropic
semi-infinite subspace $S\in Gr_{-\frac{\infty}{2}}^{(n)}$ and an element
$\omega$ from an open cone in $\hbar S/S$ we get some Frobenius manifold
structure on $\mathcal{U}$. Let us specify $S$ and $\omega$ giving rise to
initial Frobenius manifold.

We set $S:=\phi(x)H_{-}^{(n)}$, $\omega=[\phi_{1}]\in\hbar S/S$ where
$\phi_{1}=\sum_{i}\phi_{1}^{i}\rho_{i}$. It follows from definition that
subspace $S$ is opposite to $\mathcal{L}(x)$ and isotropic. Notice that
because of equation (\ref{dpsi}) $S$ and $\omega$ do not depend on
$x\in\mathcal{U}$. Also we have $\phi_{1}\in\mathcal{L}(x)\cap\hbar S$,
$i_{S}(\phi_{1})=\omega$ . Restriction of $Symbol(\frac{\partial}{\partial
x})$ on $i_{\mathcal{L}}(\phi_{1})$ is an isomorphism because $C_{\alpha
1}^{\gamma}=\delta_{\alpha}^{\gamma}$.

\begin{proposition}
Applying the construction from corollary \ref{corfrob} to the data
$\mathcal{L}(x),S,\omega$ one gets back the initial Frobenius manifold
structure on $\mathcal{U}$.
\end{proposition}

\begin{proof}
We have $\psi(x)=\phi_{1}(x)$. Notice that $[\hbar^{-1}\phi_{\alpha}]$,
$\phi_{\alpha}=\sum_{i}\phi_{\alpha}^{i}\rho_{i}$ are constant as elements of
$S/\hbar^{-1}S$ because $\frac{\partial}{\partial x^{\beta}}\phi_{\alpha}%
\in\hbar^{-1}S$. The set $\{[\hbar^{-1}\phi_{\alpha}]\}_{\alpha\in
\lbrack1,\ldots,n]\,}$ is a basis in $S/\hbar^{-1}S$. Recall that
$\frac{\partial}{\partial x^{\alpha}}\psi=\hbar^{-1}\phi_{\alpha}$. Therefore
$[\frac{\partial}{\partial x^{\alpha}}\psi]\in S/\hbar^{-1}S$ \thinspace are
constant. Therefore affine coordinates induced on $\mathcal{U}$ via the map
$[\psi(x)-\omega]:$ $\mathcal{U}\rightarrow S/\hbar^{-1}S$ coincide with
$x^{\alpha}$. Also it follows that equations (\ref{eqprfrob}) hold with
tensors $C_{\alpha\beta}^{\gamma}(x^{\alpha}),\,\eta_{\alpha\beta}%
,\frac{\partial}{\partial x^{1}}$ which coincide with tensors of the initial
Frobenius manifold.
\end{proof}

It is interesting to note that element $\psi(x)=\mathcal{L}(x)\cap(S+\omega)$
has nice meaning in all situations where Frobenius manifolds appear. That is
depending on the context it can be Baker function over small phase space of
$n$KdV hierarchy, Gromov-Witten 2-point descendent correlator, Saito primitive
form, solution to Riemann-Hilbert problem, period vector of Calabi-Yau
manifold (eventually of non-commutative deformation of it) and so on. 

\section{Lax operators.\label{slax}}

We saw above that Frobenius manifolds can be described in terms of geometry of
semi-infinite subspaces similar to the one arising in approach to integrable
hierarchies from \cite{SW}. Recall that in the latter context equations of KP
hierarchy are described in terms of flow given by multiplication by $\exp
(\sum_{r=1}^{+\infty}\hbar^{-r}x^{r})$ which acts on set of semi-infinite
subspaces. Pursuing such analogy one can conjecture that at least for certain
classes of $\,$Frobenius manifolds there exists natural \emph{enlarged} family
of semi-infinite subspaces $\mathcal{L}(x)\in$ $Gr_{\frac{\infty}{2}}^{(n)}$
depending on infinite number of parameters $\{x^{\alpha,r}\},$ $\alpha
\in\lbrack1,\ldots,n],$ $r\in\mathbb{N}$, (higher times) such that
\begin{equation}
\frac{\partial}{\partial x^{\alpha,r}}\mathcal{L}(x)\subseteq\hbar
^{-r}\mathcal{L}(x)\label{hightimes}%
\end{equation}
Intersecting $\mathcal{L}(x)$ with $\hbar S$ for some opposite subspace $S$ we
see that such family is the same as  infinite set of commuting operators
of\thinspace\ form $\frac{\partial}{\partial x^{\alpha,r}}-\sum_{j=1}^{r}%
\hbar^{-j}A_{j,(\alpha,r)}$ extending the initial set of $n$ commuting
operators $\frac{\partial}{\partial x^{\alpha}}-\hbar^{-1}C_{\alpha}$. We will
show below how dressing method can be used in order to construct and classify
such sets of operators in the case of semi-simple Frobenius manifolds.

Let $(C_{\alpha\beta}^{\gamma}(x)_{x\in\mathcal{U}},$ $\eta_{\alpha\beta},$
$e)$ define a semi-simple Frobenius manifold structure on an open domain
$\mathcal{U}\subset\mathbb{C}^{n}$. Recall that unless the converse is
explicitly mentioned we do not assume it to be conformal, i.e. we do not
assume existence of an Euler vector field. Three-tensor $C_{\alpha\beta
}^{\gamma}(x)$ defines structure constants of commutative associative
multiplication on tangent space $T_{x}\mathcal{U}$ for any $x\in\mathcal{U}$ .
We denote this multiplication via $"\circ"$. Multiplication operators
$"v\circ"\in End(T_{x}\mathcal{U}),$ $v\in T_{x}\mathcal{U}$ form commutative
$n-$dimensional subalgebra in $End(T_{x}\mathcal{U})$. Our basic assumption
(semisimplicity) is that for any $x\in\mathcal{U}$ we have decomposition
$T_{x}\mathcal{U}=\oplus_{i=1}^{n}\theta_{i}$ where $\theta_{i}$ are
eigenspaces of all operators $"v\circ"$, $v\in T_{x}\mathcal{U}$ simultaneously.

Let us introduce our Lax operator. Contrary to the usual case it is in fact a
set of $n$ commuting operators in $n$ variables:
\begin{equation}
L=\{L_{\alpha}(x)\}_{\alpha\in\lbrack1,\ldots,n]}\,\,\,\,\,L_{\alpha}%
=\frac{\partial}{\partial x^{\alpha}}-\hbar^{-1}C_{\alpha}(x)\label{Lalpha}%
\end{equation}
$[L_{\alpha},L_{\beta}]=0$, where $C_{\alpha}$ are functions of $x=\{x^{\alpha
}\}_{\alpha\in\lbrack1,\ldots,n]}$ with values in $Mat(n,\mathbb{C)}$ whose
entries are $(C_{\alpha})_{\beta}^{\gamma}=C_{\alpha\beta}^{\gamma}(x)$. We
first apply some transformations to reduce the set of operators $\{L_{\alpha
}(x)\}$ to a simpler form. Although it might look strange in the view of
standard definition of Frobenius manifold, it is natural from the point of
view described in previous section  to consider separately changes of
coordinates on $\mathcal{U}$ acting on $\{L_{\alpha}(x)\}$ as on components of
$Mat(n,\mathbb{C)-}$connection in \textit{trivial }bundle (and not as on
components of a connection on $T\mathcal{U}$ ) and gauge transformations
changing the choice of frame in the trivial bundle.

Let us choose nonzero vectors $v_{i}\in\theta_{i}$ where $\theta_{i}\subset
T_{x}\mathcal{U}$ are eigenspaces of multiplication operators, $v_{i}=\sum
v_{i}^{\beta}(x)\frac{\partial}{\partial x^{\beta}}$ and put $(T_{0}%
(x))_{i}^{\beta}=$ $v_{i}^{\beta}(x)$, $T_{0}(x)\in GL(n,\mathbb{C})$. Then
for all $\alpha$ simultaneously we have: $T_{0}^{-1}C_{\alpha}T_{0}=a_{\alpha
}(x)$ where $a_{\alpha}(x)=diag(a_{\alpha}^{1},\ldots,a_{\alpha}^{n})$ is a
function of $x$ with values in $Diag$. Note that linear combinations of
$a_{\alpha}(x)$ span $Diag$ as linear space for any $x\in\mathcal{U}$.
Transformation $T_{0}(x)$ is determined uniquely up to right multiplication by
function with values in $Diag$. Then transformed Lax operator $\widetilde
{L}=\{\widetilde{L}_{\alpha}\},\widetilde{L}_{\alpha}:=T_{0}^{-1}L_{\alpha
}T_{0}$ has the form
\begin{equation}
\widetilde{L}_{\alpha}=\frac{\partial}{\partial x^{\alpha}}+q_{\alpha}%
-\hbar^{-1}a_{\alpha}(x)\label{Ltildealpha}%
\end{equation}
with $q_{\alpha}=T_{0}^{-1}(\frac{\partial}{\partial x^{\alpha}}T_{0})$.
Vectors $v_{i}$ satisfy $\eta(v_{i},v_{j})=0$ for $i\neq j$ where $\eta$ is
the pairing defined by tensor $\eta_{\alpha\beta}$. We can normalize $v_{i}$
so that $\eta(v_{i},v_{i})=1$. Then $T_{0}T_{0}^{\top}=Id$ and $(T_{0}%
^{-1}(\frac{\partial}{\partial x^{\alpha}}T_{0}))^{\top}=-T_{0}^{-1}%
(\frac{\partial}{\partial x^{\alpha}}T_{0})$. We denote transformed Lax
operator written using such orthonormal normalization $\widetilde{L}%
^{orth}=\{\widetilde{L}_{\alpha}^{orth}\}$

We can simplify further the form of our Lax operators by appropriate change of
coordinates on $\mathcal{U}$. The equation $[\widetilde{L}_{\alpha}%
,\widetilde{L}_{\beta}]=0$ implies that $\frac{\partial}{\partial x^{\beta}%
}a_{\alpha}(x)-\frac{\partial}{\partial x^{\alpha}}a_{\beta}(x)=[T_{0}%
^{-1}(\frac{\partial}{\partial x^{\beta}}T_{0}),a_{\alpha}]+[T_{0}^{-1}%
(\frac{\partial}{\partial x^{\alpha}}T_{0}),a_{\beta}]$. Notice that diagonal
entries of $[B,a]$ are zero for $a\in Diag$. Therefore $\frac{\partial
}{\partial x^{\beta}}a_{\alpha}(x)-\frac{\partial}{\partial x^{\alpha}%
}a_{\beta}(x)=0$ and $a_{\alpha}(x)=\frac{\partial}{\partial x^{\alpha}}u(x)$
for some $u(x)=diag(u^{1},\ldots,u^{n})$ which is a function of $x^{\alpha}$
with values in $Diag$. Notice that $u^{1},\ldots,u^{n}$ define a new set of
coordinates on $\mathcal{U}$ since $\det\frac{\partial u^{i}(x)}{\partial
x^{\alpha}}=\det a_{\alpha}^{i}\neq0$. In this coordinates we have
\begin{equation}
\widetilde{L}=\{\widetilde{L}_{i}(u)\},\widetilde{L}_{i}=\sum_{\alpha}%
\frac{\partial x^{\alpha}}{\partial u^{i}}\widetilde{L}_{\alpha}%
=\frac{\partial}{\partial u^{i}}+q_{i}(u)-\hbar^{-1}e_{i}\label{Li}%
\end{equation}
where $e_{i}$ is the constant matrix whose entries are $(e_{i})_{jk}%
=\delta_{ij}\delta_{ik}$. We have $[q_{i}(u),e_{j}]-[q_{j}(u),e_{i}]=0$ as a
consequence of $[\widetilde{L}_{i},\widetilde{L}_{j}]=0.$ It follows that
$q_{i}(x)$ may have non-zero entries on $i-$th row, $i$-th column and diagonal
only and that for some matrix $Q_{kl}(u)$ we have $q_{i}=[e_{i},Q]+\sum
_{k}(q_{i})_{kk}e_{k}$ simultaneously for all $q_{i}$. Matrix $Q_{kl}$ is
determined uniquely if one requires that its diagonal entries are zero. For
$\widetilde{L}^{orth}$ we have $(q_{i}(u))^{\top}=-q_{i}(u)$ and
$q_{i}(u)=[e_{i},Q]$ for  $Q_{kl}(u),Q^{\top}=Q,$ $Q_{kk}=0$.

Notice that in order to reduce commuting operators of form $\frac{\partial
}{\partial x^{\alpha}}-\hbar^{-1}C_{\alpha}(x)$ to operators of form
(\ref{Ltildealpha}) and (\ref{Li}) we have used only that $C_{\alpha}(x)$ are
simultaneously diagonalizable and that the vector space spanned by their
spectrums coincides with the space of all diagonal matrices at any
$x\in\mathcal{U}$.

\section{Dressing transformations.}

\begin{proposition}
\label{dress}There exists formal power series $\widetilde{T}(\hbar
,u)=Id+\sum_{k=1}^{\infty}\hbar^{k}\widetilde{T}_{k}(u)$ such that for all
$i\in\lbrack1,\ldots,n]$ \emph{simultaneously}
\begin{equation}
\widetilde{T}^{-1}\widetilde{L}_{i}\widetilde{T}=\frac{\partial}{\partial
u^{i}}-e_{i}\hbar^{-1}+\sum_{k=0}^{\infty}\hbar^{k}h_{k,i}(u)\label{Ttilde}%
\end{equation}
where $h_{k,i}(u)$ are functions with values in $Diag$. $\widetilde{T}$ is
determined up to right multiplication by function with values in
$Id+(Diag)\hbar\lbrack\lbrack\hbar]]$. If $\widetilde{T}$ $\,$is represented
\thinspace in the form $\widetilde{T}=(id+\hbar S_{1}(u))\circ\ldots
\circ(id+\hbar^{k}S_{k}(u))\circ\ldots$ then $\widetilde{T}$ is determined
uniquely by requiring that diagonal entries of $S_{k}$ are zero. In this case
entries of $S_{k}$ (and consequently of $\widetilde{T}_{k}$) can be found
recursively and they are given by certain differential polynomials in
$(q_{i})_{kl}$ with zero free terms. For $\widetilde{L}^{orth}$ we have
$\widetilde{T}(\hbar)\widetilde{T}(-\hbar)^{\top}=Id$ , $h_{k,2i}(u)=0$.
\end{proposition}

\begin{proof}
We look for $\widetilde{T}$ in the form $\widetilde{T}=(id+\hbar
S_{1}(u))\circ\ldots\circ(id+\hbar^{k}S_{k}(u))\circ\ldots$ so that
\begin{multline*}
\lbrack(id+\hbar S_{1}(u))\circ\ldots\circ(id+\hbar^{k}S_{k}(u))]^{-1}%
\widetilde{L}_{i}[(id+\hbar S_{1}(u))\circ\ldots\circ(id+\hbar^{k}%
S_{k}(u))]=\\
=\frac{\partial}{\partial u^{i}}-e_{i}\hbar^{-1}+\sum_{j=0}^{k-1}\hbar
^{j}h_{j,i}(u)+\sum_{j=k}^{\infty}\hbar^{j}H_{j,i}(u)
\end{multline*}
where $h_{j,i}(u)\in Diag$. Then equations for $S_{k+1}(u),$ $h_{k,i}(u)$ are
\[
-e_{i}\circ S_{k+1}+H_{k,i}=-S_{k+1}\circ e_{i}+h_{k,i}%
\]
Notice that $[\widetilde{L}_{i},\widetilde{L}_{j}]=0$ implies same relation
for Lax operators conjugated by $(id+\hbar S_{1}(u))\circ\ldots\circ
(id+\hbar^{k}S_{k}(u))$. Therefore $[e_{i},H_{k,j}]-[e_{j},H_{k,i}]\in Diag$.
It follows that $[e_{i},H_{k,j}]-[e_{j},H_{k,i}]=0$. Therefore matrix
$H_{k,i}$ may have non-zero entries on $i-$th row, $i$-th column and diagonal
only and for some matrix $S_{k+1}$  and some matrices $h_{k,i}\in Diag$ we
have $H_{k,i}=[e_{i},S_{k+1}]+h_{k,i}$. Such matrix $S_{k+1}$  is determined
uniquely by requirement that its diagonal entries are zero. It is easy to see
by induction that entries of $H_{k,i}$ and $S_{k}$ are then differential
polynomials in $q_{j}$ with zero free terms and therefore the same holds for
$\widetilde{T}_{k}(u)$. Let $\widetilde{T}^{\prime}$ is another transformation
such that $\widetilde{T}^{\prime-1}\widetilde{L}_{i}\widetilde{T}^{\prime
}=\frac{\partial}{\partial u^{i}}-e_{i}\hbar^{-1}+h_{i}^{\prime},$
$h_{i}^{\prime}\in Diag[[\hbar]]$. Let us set $R=\widetilde{T}^{-1}%
\widetilde{T}^{\prime}$, $R=Id+\sum_{k=1}^{\infty}\hbar^{k}R_{k}$ then
$R^{-1}(\frac{\partial}{\partial u^{i}}-e_{i}\hbar^{-1}+h_{i})R=\frac
{\partial}{\partial u^{i}}-e_{i}\hbar^{-1}+h_{i}^{\prime}$, $h_{i},$
$h_{i}^{\prime}\in Diag[[\hbar]]$. Let $l$ be such that $R_{k}\in Diag$ for
all $k<l$ then $[e_{i},R_{l}]\in Diag$, therefore $[e_{i},R_{l}]=0$ and
consequently $R_{l}\in Diag$. Notice that for $\widetilde{L}^{orth}$ we have
by induction $(H_{k,i})^{\top}=(-1)^{k+1}H_{k,i}$ and therefore $(S_{k+1}%
)^{\top}=(-1)^{k}S_{k+1}$and $h_{k,2i}(u)=0$.
\end{proof}

\begin{proposition}
\label{cordress}\bigskip If one sets $T(\hbar,x)=T_{0}(x)\circ\widetilde
{T}(u(x))$ then for all $\alpha\in\lbrack1,\ldots,n]$ \emph{simultaneously}
\begin{equation}
T^{-1}L_{\alpha}T=\widetilde{T}^{-1}(u(x))\widetilde{L}_{\alpha}\widetilde
{T}(u(x))=\frac{\partial}{\partial x^{\alpha}}-a_{\alpha}(x)\hbar^{-1}%
+\sum_{k=0}^{+\infty}\hbar^{k}h_{k,\alpha}(x)\label{T}%
\end{equation}
where $h_{k,\alpha}(x):=\sum_{i}a_{\alpha}^{i}h_{k,i}(u(x))$ and $a_{\alpha
}(x)$ are functions with values in $Diag$. Transformation $T(x)\in
Mat(n,\mathbb{C})[[\hbar]]$ satisfying $T^{-1}L_{\alpha}T=\frac{\partial
}{\partial x^{\alpha}}+h_{\alpha}$, $h_{\alpha}\in\hbar^{-1}Diag[[\hbar]]$ for
all $\alpha\in\lbrack1,\ldots,n]$ simultaneously\emph{\ }is defined uniquely
up to right multiplication by a function with values in $Diag[[\hbar]]$.
\end{proposition}

\begin{proof}
The first statement is a consequence of the previous proposition. The proof of
the second statement is the same as proof of analogous statement from the
previous proposition.
\end{proof}

Below we set $\widetilde{L}_{i}^{norm}:=\widetilde{T}^{-1}\widetilde{L}%
_{i}\widetilde{T}$ , $L_{\alpha}^{norm}:=T^{-1}L_{\alpha}T$. We also set
$h_{-1,\alpha}=-a_{\alpha}$. We adopt terminology which calls such
transformations reducing Lax operators to diagonal form dressing transformations.

\section{Lax equations.}

The construction proceeds quite analogously to the standard case of Lax
operator in one variable (see \cite{DS}, \S1).

Notice that if $b\in Diag((\hbar))$ then $[\widetilde{L}_{i}^{norm},b]=0$,
$[L_{\alpha}^{norm},b]=0$. Therefore $[\widetilde{L}_{i},\widetilde
{T}b\widetilde{T}^{-1}]=0$, $[L_{\alpha},TbT^{-1}]=0$. For any $b\in
Diag((\hbar))$ we set $\widetilde{\varphi}(b)=$ $\widetilde{T}b\widetilde
{T}^{-1}$, $\varphi(b)=TbT^{-1}$.

\begin{proposition}
For any $b\in Diag((\hbar))$ we have
\begin{align*}
\lbrack L_{\alpha},\varphi(b)_{\leq-1}]  &  =\hbar^{-1}[C_{\alpha}%
,\varphi(b)_{0}]\\
\lbrack\widetilde{L}_{i},\widetilde{\varphi}(b)_{\leq0}]  &  =\hbar^{0}%
[e_{i},\widetilde{\varphi}(b)_{1}]
\end{align*}
\end{proposition}

\begin{proof}
$[\widetilde{L}_{i},\widetilde{\varphi}(b)_{\leq0}]\in Mat(n,\mathbb{C}%
)[\hbar^{-1}]$, on the other hand $[\widetilde{L}_{i},\widetilde{\varphi
}(b)_{\leq0}]=-$ $[\widetilde{L}_{i},\widetilde{\varphi}(b)_{>0}]=\hbar
^{0}[e_{i},\widetilde{\varphi}(b)_{1}]+s$, $s\in Mat(n,\mathbb{C})\hbar
\lbrack\lbrack\hbar]]$. It follows that $s=0$. The proof in the case of
$L_{\alpha}$ is the same.
\end{proof}

\bigskip Our basic set of equations is
\begin{equation}
\frac{\partial L_{\alpha}}{\partial t}=[\sum_{j=1}^{m}\hbar^{-j}%
A_{j},L_{\alpha}]\label{Lalpha(t)}%
\end{equation}
where $L_{\alpha}=\frac{\partial}{\partial x^{\alpha}}-\hbar^{-1}C_{\alpha}$,
$\alpha\in\lbrack1,\ldots,n]$, $[L_{\alpha},L_{\beta}]=0$, where $C_{\alpha}$
and $A_{j}$ are functions of $x^{\alpha},t$ with values in $Mat(n,\mathbb{C)}$
and $\sum_{j=1}^{m}\hbar^{-j}A_{j}=\varphi(b)_{\leq-1}$ for some fixed $b\in
Diag((\hbar))$. Recall that $\varphi(b):=T(\hbar,x,t)bT^{-1}(\hbar,x,t)$ where
$T$ is dressing transformation from proposition \ref{cordress} (with $t$
considered as a parameter). Notice that in spite of non-uniqueness in the
choice of $T$ (recall that it is defined up to right multiplication by $R\in
Diag[[\hbar]]$) $\varphi(b)$ is well-defined. In fact $\varphi(b)$ for
$b=\sum_{j=-\infty}^{j=m}$ $\hbar^{-j}b_{j}$ does not depend on $b_{j}$ with
$j\leq0$. Notice also that $[[\sum_{j=1}^{m}\hbar^{-j}A_{j},L_{\alpha
}],L_{\beta}]+[L_{\alpha},[\sum_{j=1}^{m}\hbar^{-j}A_{j},L_{\beta}]]=0$ which
implies that it is enough to impose condition $[L_{\alpha},L_{\beta}]=0$ at
some initial value of $t$ only.

We also consider analogous equations for $\{\widetilde{L}_{i}\}$ $:$%
\begin{equation}
\frac{\partial\widetilde{L}_{i}}{\partial t}=[\sum_{j=0}^{m}\hbar
^{-j}\widetilde{A}_{j},\widetilde{L}_{i}]\label{Li(t)}%
\end{equation}
where $\widetilde{L}_{i}=\frac{\partial}{\partial u^{i}}+q_{i}(u,t)-\hbar
^{-1}e_{i}$, $i=1,\ldots,n$, $[\widetilde{L}_{i},\widetilde{L}_{j}]=0$,
$q_{i}(u,t)$ and $\widetilde{A}_{j}(u,t)$ are functions of $u^{i}$, $t$ with
values in $Mat(n,\mathbb{C)}$; $\sum_{j=0}^{m}\hbar^{-j}\widetilde{A}%
_{j}(u,t)=\widetilde{\varphi}(b)_{\leq0}$ for some fixed $b\in Diag((\hbar))$,
$\widetilde{\varphi}(b)=\widetilde{T}(\hbar,u,t)b\widetilde{T}^{-1}%
(\hbar,u,t)$ where $\widetilde{T}(\hbar,u,t)$ is dressing transformation for
$\widetilde{L}_{i}$ defined by equation (\ref{Ttilde}) with $t$ as a
parameter. Also $\widetilde{\varphi}(b)$ for  $b=\sum_{j=-\infty}^{j=m}$
$\hbar^{-j}b_{j}$ depends only on $b_{j}$ with $j\geq0$.

Similarly we have for $\{\widetilde{L}_{\alpha}\}$ the following equations:%

\begin{equation}
\frac{\partial\widetilde{L}_{\alpha}}{\partial t}=[\sum_{j=0}^{m}\hbar
^{-j}\widetilde{A}_{j},\widetilde{L}_{\alpha}]\label{Ltildealpha(t)}%
\end{equation}
where $\widetilde{L}_{\alpha}=\frac{\partial}{\partial x^{\alpha}}+q_{\alpha
}(x,t)-\hbar^{-1}a_{\alpha}(x)$, $\alpha\in\lbrack1,\ldots,n]$, $[\widetilde
{L}_{\alpha},\widetilde{L}_{\beta}]=0$, $a_{\alpha}(x)\in Diag$ , $q_{\alpha
}(x,t)$ and $\widetilde{A}_{j}(x,t)$ are functions of $x^{\alpha}$, $t$ with
values in $Mat(n,\mathbb{C)}$, $\sum_{j=0}^{m}\hbar^{-j}\widetilde{A}%
_{j}(x,t)=\widetilde{\varphi}(b)_{\leq0}$ for some fixed $b\in Diag((\hbar))$,
$\widetilde{\varphi}(b)=\widetilde{T}(\hbar,x,t)b\widetilde{T}^{-1}%
(\hbar,x,t)$ and $\widetilde{T}(\hbar,x,t)$ is dressing transformation defined
by equation (\ref{T}) with $t$ as a parameter. As above $\widetilde{\varphi
}(b)$ for $b=\sum_{j=-\infty}^{j=m}$ $\hbar^{-j}b_{j}$ depends only on $b_{j}$
with $j\geq0$.

These three sets of equations are in fact closely related.

\begin{proposition}
Let operators $\{\widetilde{L}_{\alpha}(x,t_{1})\}$ of form (\ref{Ltildealpha}%
) are related with operators $\{\widetilde{L}_{i}(u,t_{1})\}$ of form
(\ref{Li}) by a change of coordinates $u^{i}=u^{i}(x^{\alpha})$ so that
$\widetilde{L}_{\alpha}(x,t_{1})=\sum_{i}\frac{\partial u^{i}}{\partial
x^{\alpha}}\widetilde{L}_{i}(u(x),t_{1})$, then $\{\widetilde{L}_{i}(u,t)\}$
satisfies equations (\ref{Li(t)}) for some $b\in Diag[\hbar^{-1}]$ iff
$\{\widetilde{L}_{\alpha}(x,t)=\sum_{i}\frac{\partial u^{i}}{\partial
x^{\alpha}}\widetilde{L}_{i}(u(x),t)\}$ satisfies equations (\ref{Lalpha(t)})
with the same $b\in Diag[\hbar^{-1}]$.
\end{proposition}

\begin{proof}
If operators $\widetilde{L}_{i}(u,t)$ satisfy equations (\ref{Li(t)}) then
their linear combinations $\sum_{i}\frac{\partial u^{i}}{\partial x^{\alpha}%
}\widetilde{L}_{i}(u,t)$ with coefficients which are independent of $t$
satisfy
\[
\frac{\partial\widetilde{L}_{\alpha}(x,t)}{\partial t}=\sum_{i}\frac{\partial
u^{i}}{\partial x^{\alpha}}\frac{\widetilde{L}_{i}(u(x),t)}{\partial t}%
=\sum_{i}[\widetilde{\varphi}(b)_{\leq0},\frac{\partial u^{i}}{\partial
x^{\alpha}}\widetilde{L}_{i}(u(x),t)]
\]
Also if $\widetilde{T}(\hbar,u,t)$ is dressing transformation for
$\widetilde{L}_{i}(u,t)$: $\widetilde{T}^{-1}\widetilde{L}_{i}\widetilde
{T}=\frac{\partial}{\partial u^{i}}-e_{i}\hbar^{-1}+\,h_{i}$, $h_{i}\in
Diag[[\hbar]]$, then $\widetilde{T}(\hbar,x,t)=\widetilde{T}(\hbar,u(x),t)$ is
dressing transformation for $\widetilde{L}_{\alpha}(x,t)$:
\[
\widetilde{T}^{-1}[\sum_{i}\frac{\partial u^{i}}{\partial x^{\alpha}%
}\widetilde{L}_{i}(u(x),t)]\widetilde{T}=\frac{\partial}{\partial x^{\alpha}%
}-a_{\alpha}(x)\hbar^{-1}+h_{\alpha}(x,t),\,h_{\alpha}(x,t)\in Diag[[\hbar]]
\]
Conversely if $\widetilde{L}_{\alpha}(x,t)$ satisfy equations (\ref{Lalpha(t)}%
) then $[\widetilde{\varphi}(b)_{\leq0},\widetilde{L}_{\alpha}%
(x,t)]=[\widetilde{\varphi}(b),a_{\alpha}]$ and $\frac{\partial}{\partial
t}a_{\alpha}=0$. Thus $\{\widetilde{L}_{i}(u,t):=\sum_{\alpha}\frac{\partial
x^{\alpha}}{\partial u^{i}}\widetilde{L}_{\alpha}(x(u),t)\}$ has the form
(\ref{Li}) for any $t$ and satisfies equations.(\ref{Li(t)}).
\end{proof}

We saw above that given operators $\{L_{\alpha}(x)\}$ in the form
(\ref{Lalpha}) there exists a gauge transformation $T_{0}(x)$ such that
operators $\widetilde{L}_{\alpha}(x)=T_{0}^{-1}L_{\alpha}(x)T_{0}$ are in the
form (\ref{Ltildealpha}). Conversely, given operators $\{\widetilde{L}%
_{\alpha}(x)\}$ of form (\ref{Ltildealpha}) we have $[\frac{\partial}{\partial
x^{\alpha}}+q_{\alpha},\frac{\partial}{\partial x^{\beta}}+q_{\beta}]=0$.
Therefore there exists a gauge transformation $S_{0}(x)$ such that
$S_{0}(\frac{\partial}{\partial x^{\alpha}}+q_{\alpha})S_{0}^{-1}%
=\frac{\partial}{\partial x^{\alpha}}$ and consequently $S_{0}\widetilde
{L}_{\alpha}(x)S_{0}^{-1}$ is in the form (\ref{Lalpha}). 

\begin{proposition}
Let operators $\{L_{\alpha}(x)\}_{\alpha\in\lbrack1,\ldots,n]}$ of form
(\ref{Lalpha}) are related with operators $\{\widetilde{L}_{\alpha}(x)\}$ of
form (\ref{Ltildealpha}) by a gauge transformation: $\widetilde{L}_{\alpha
}(x)=T_{0}^{-1}(x)L_{\alpha}(x)T_{0}(x)$. Let $\{L_{\alpha}(x,t)\}_{\alpha
\in\lbrack1,\ldots,n]},$ $L_{\alpha}(x,t_{1})=L_{\alpha}(x)$, satisfy
equations (\ref{Lalpha(t)}) for some $b\in Diag[\hbar^{-1}]$. If one sets
$R(x,t)$ to be family of gauge transformations $R(x,t)$ such that $\partial
R/\partial t=-\varphi(b)_{0}R$, $R(x,t_{1})=T_{0}(x)$, then $\widetilde
{L}_{\alpha}(x,t):=R^{-1}(x,t)L_{\alpha}(x,t)R(x,t)$ have the form
(\ref{Ltildealpha}) for any $t$ and satisfy equations (\ref{Ltildealpha(t)})
with the same $b\in Diag[\hbar^{-1}]$. Conversely if $\{\widetilde{L}_{\alpha
}(x,t)\}_{\alpha\in\lbrack1,\ldots,n]}$, $\widetilde{L}_{\alpha}%
(x,t_{1})=\widetilde{L}_{\alpha}(x)$, satisfy equations.(\ref{Ltildealpha(t)})
for some $b\in Diag[\hbar^{-1}]$, then for $R(x,t)$ such that $\partial
R/\partial t=-R[\widetilde{\varphi}(b)]_{0}$, $R(x,t_{1})=T_{0}(x)$, operators
$L_{\alpha}(x,t):=R(x,t)\widetilde{L}_{\alpha}(x,t)R^{-1}(x,t)$ are of form
(\ref{Lalpha}) and satisfy equations (\ref{Lalpha(t)}) with the same $b\in
Diag[\hbar^{-1}]$.
\end{proposition}

\begin{proof}
We have $T_{0}^{-1}C_{\alpha}(x,t_{1})T_{0}=a_{\alpha}(x)$, $q_{\alpha}%
(t_{1})=T_{0}^{-1}(\frac{\partial}{\partial x^{\alpha}}T_{0})$. If $\partial
L_{\alpha}/\partial t=[\varphi(b)_{\leq-1},L_{\alpha}]$ then
\[
\frac{\partial C_{\alpha}}{\partial t}=-[\varphi(b)_{0},C_{\alpha}]
\]
and therefore
\begin{multline*}
\frac{\partial(R^{-1}C_{\alpha}R)}{\partial t}=-R^{-1}\frac{\partial
R}{\partial t}R^{-1}C_{\alpha}R+R^{-1}\frac{\partial C_{\alpha}}{\partial
t}R+R^{-1}C_{\alpha}\frac{\partial R}{\partial t}=\\
=R^{-1}\varphi(b)_{0}C_{\alpha}R-R^{-1}\varphi(b)_{0}C_{\alpha}R+R^{-1}%
C_{\alpha}\varphi(b)_{0}R-R^{-1}C_{\alpha}\varphi(b)_{0}R=0
\end{multline*}
Therefore $R^{-1}C_{\alpha}(x,t)R=a_{\alpha}(x)$ does not depend on $t$ and
$R^{-1}L_{\alpha}(x,t)R=\frac{\partial}{\partial x^{\alpha}}+q_{\alpha
}(x,t)-\hbar^{-1}a_{\alpha}(x)$, $q_{\alpha}=R^{-1}\frac{\partial R}{\partial
x^{\alpha}}$, $\alpha\in\lbrack1,\ldots,n]$, are of the form
(\ref{Ltildealpha}) for any $t$. If $\widetilde{T}(x,t)$ is the dressing
transformation for $\frac{\partial}{\partial x^{\alpha}}+q_{\alpha}%
(x,t)-\hbar^{-1}a_{\alpha}(x)$ then $R\circ\widetilde{T}$ is dressing
transformation for $L_{\alpha}(x,t)$ and $R^{-1}\varphi(b)R=$ $=\widetilde
{T}b\widetilde{T}^{-1}=\widetilde{\varphi}(b)$.We have
\[
R^{-1}\circ(\frac{\partial}{\partial t})\circ R=\frac{\partial}{\partial
t}+R^{-1}\frac{\partial R}{\partial t}=\frac{\partial}{\partial t}%
-\widetilde{\varphi}(b)_{0}%
\]

and%

\[
0=[R^{-1}\circ(\frac{\partial}{\partial t})\circ R,R^{-1}\circ(\frac{\partial
}{\partial x^{\alpha}})\circ R]=[\frac{\partial}{\partial t}-\widetilde
{\varphi}(b)_{0},\frac{\partial}{\partial x^{\alpha}}+q_{\alpha}(x,t)]
\]
Therefore $\partial q_{\alpha}/\partial t=[\widetilde{\varphi}(b)_{0}%
,\frac{\partial}{\partial x^{\alpha}}+q_{\alpha}(x,t)]=[\widetilde{\varphi
}(b)_{\leq0},\frac{\partial}{\partial x^{\alpha}}+q_{\alpha}(x,t)-\hbar
^{-1}a_{\alpha}(x)]$. The proof of converse statement uses similar arguments.
\end{proof}

Equations (\ref{Lalpha(t)}) describe essentially the largest possible
extension of our family of semi-infinite subspaces to the family satisfying
(\ref{hightimes}) because of the following proposition.

\begin{proposition}
Let $[\sum_{j=1}^{m}\hbar^{-j}A_{j}(x),L_{\alpha}]\in\hbar^{-1}%
Mat(n,\mathbb{C})$ for all $\alpha\in\lbrack1,\ldots,n]$, then $\sum_{j=1}%
^{m}\hbar^{-j}A_{j}(x)=\varphi(b)_{\leq-1}$ for some $b=\sum_{j=1}^{m}%
b_{j}\hbar^{-j}\in Diag[\hbar^{-1}]$ and for $j>1$: $b_{j}=const$.
\end{proposition}

\begin{proof}
Let us set $T(\sum_{j=1}^{m}\hbar^{-j}A_{j}(x))T^{-1}=\sum_{j=-\infty}%
^{m}\hbar^{-j}B_{j}$. We have
\begin{equation}
\lbrack\sum_{j=-\infty}^{m}\hbar^{-j}B_{j},\frac{\partial}{\partial x^{\alpha
}}-\hbar^{-1}a_{\alpha}+h_{\alpha}]\in\hbar^{-1}Mat(n,\mathbb{C}%
),\,\,\text{where }h_{\alpha}\in Diag[[\hbar]]\label{eqfull}%
\end{equation}
Let $B_{l}\in Diag$ for all $l>j>0$. Then $[B_{j},a_{\alpha}]\in Diag$ and
therefore $[B_{j},a_{\alpha}]=0.$ It follows that $B_{j}\in Diag$ and the same
is true for all $B_{j}$ with $j>0$. Therefore using again (\ref{eqfull})
$\frac{\partial}{\partial x^{\alpha}}B_{j}=0$ for all $j>1$.
\end{proof}

Below we will mainly consider equations (\ref{Lalpha(t)}).

\subsection{Integrals of motion.}

Let $L_{\alpha}(x,t)$ satisfy (\ref{Lalpha(t)}) and $T^{-1}(x,t)L_{\alpha
}T(x,t)=\frac{\partial}{\partial x^{\alpha}}+\sum_{k=-1}^{+\infty}\,\hbar
^{-k}h_{k,\alpha}(x,t)$, $h_{k,\alpha}(x,t)\in Diag$.

\begin{proposition}
\label{integrals} $\,h_{k,\alpha}(x,t)$ satisfy $\,\partial h_{k,\alpha
}(x,t)/\partial t+\partial B_{k}/\partial x^{\alpha}$ for some $B_{k}(x,t)\in
Diag$.
\end{proposition}

\begin{proof}
Equations \ref{Lalpha(t)} can be written as $[\frac{\partial}{\partial t}%
-\sum_{k=1}^{k=m}\hbar^{-k}A_{k},L_{\alpha}]=0,\,\,\alpha\in\lbrack
1,\ldots,n]$. It follows that
\begin{equation}
\lbrack T^{-1}(\frac{\partial}{\partial t}-\sum_{j=1}^{j=m}\hbar^{-j}%
A_{j})T,\frac{\partial}{\partial x^{\alpha}}+\sum_{k=-1}^{+\infty}\,\hbar
^{-k}h_{k,\alpha}(x,t)]=0\label{cons}%
\end{equation}
Let $T^{-1}(\frac{\partial}{\partial t}-\sum_{j=0}^{j=m}\hbar^{-j}%
A_{j})T=\frac{\partial}{\partial t}-\sum_{j=-\infty}^{j=m}\hbar^{-j}B_{j}$.
Notice that $\{a_{\alpha}=h_{-1,\alpha}\}_{\,\,\alpha\in\lbrack1,\ldots,n]}$
generate $Diag$. Let for all $j>l$ $B_{j}\in Diag.$ Then $[B_{l},a_{\alpha
}]\in Diag$ and therefore $[B_{l},a_{\alpha}]=0$ for all $\alpha\in
\lbrack1,\ldots,n]$ and consequently $B_{l}\in Diag$. We see that $B_{j}\in
Diag$. Therefore setting $B:=\sum_{j=-\infty}^{j=m}\hbar^{-j}B_{j}$ equation
(\ref{cons}) can be written
\[
\frac{\partial h_{\alpha}}{\partial t}+\frac{\partial B}{\partial x^{\alpha}%
}=0
\]
\end{proof}

Let $x^{\alpha}=x^{\alpha}(s)$ be a curve on $\mathcal{U}$ and let us consider
evolution of $L|_{x(s)}=L(s,t)$ which is restriction on $x^{\alpha}(s)$ of the
flat connection corresponding to $\{L_{\alpha}(t)\}_{\alpha\in\lbrack
1,\ldots,n]}$, $L(s,t):=\sum_{\alpha}(\partial x^{\alpha}/\partial
s)L_{\alpha}(x,t)$. If $\{L_{\alpha}(x,t)\}_{\alpha\in\lbrack1,\ldots,n]}$
satisfy equation(\ref{Lalpha(t)}) for some $b\in Diag[\hbar^{-1}]$, then we
have analogous equation for $L(s,t)$ : $\partial L/\partial t=[\varphi
(b)_{\leq-1},L]$. Proposition \ref{integrals} implies that $\sum_{\alpha
}h_{k,\alpha}^{j}(x,t)dx^{\alpha}|_{x(s)}$ are densities of conservation laws
for this equation. Notice also that $\sum_{\alpha}h_{k,\alpha}^{j}%
(x,t)dx^{\alpha}$ is a closed form because of flatness of $L$. Therefore for a
closed curve $x(s)$ the integrals $\int_{x(s)}\sum_{\alpha}h_{k,\alpha}%
^{j}(x,t)dx^{\alpha}$ do not change under deformations of $x(s)$. The same is
true for a curve satisfying periodic or appropriate asymptotic boundary
conditions. Analogous results hold for equations (\ref{Li(t)}),
(\ref{Lalpha(t)}).

\subsection{Commutativity of flows.}

Let us consider two sets of equations:
\begin{equation}
\,\,\frac{\partial L_{\alpha}}{\partial t^{i}}=[M_{\leq-1}^{(i)},L_{\alpha
}],\,\,M^{(i)}=Tb^{(i)}T^{-1},\,\,b^{(i)}\in Diag((\hbar
)),\,i=1,2\,\,\label{dt1dt2}%
\end{equation}

\begin{proposition}
\label{pr:comm}$\,$Flows defined by $b^{(i)}\in Diag((\hbar)),\,i=1,2\,$%
\thinspace\thinspace in (\ref{dt1dt2}) commute.
\end{proposition}

\begin{proof}
We must prove that $\frac{\partial^{2}L_{\alpha}}{\partial t^{1}\partial
t^{2}}=\,\frac{\partial^{2}L_{\alpha}}{\partial t^{2}\partial t^{1}}$ where
derivatives are computed by means of equation (\ref{dt1dt2}). The proof is
parallel to the standard case of Lax operator in one variable (proposition 1.7
from \cite{DS}). We have$\frac{\partial}{\partial t^{1}}$ $\frac{\partial
L_{\alpha}}{\partial t^{2}}=[\frac{\partial}{\partial t^{1}}(M_{\leq-1}%
^{(2)}),L_{\alpha}]+[M_{\leq-1}^{(2)},[M_{\leq-1}^{(1)},L_{\alpha}]]$. By the
same arguments as in proof of proposition \ref{integrals} we have
$T^{-1}(\frac{\partial}{\partial t^{i}}-M_{\leq-1}^{(i)})T\in Diag((\hbar))$.
Therefore $[T^{-1}(\frac{\partial}{\partial t^{i}}-M_{\leq-1}^{(i)}%
)T,b^{(j)}]=0$ and $[\frac{\partial}{\partial t^{i}}-M_{\leq-1}^{(i)}%
,M^{(j)}]=0$. Therefore $\frac{\partial}{\partial t^{i}}M_{\leq-1}%
^{(j)}=[M_{\leq-1}^{(i)},M^{(j)}]_{\leq-1}$and
\begin{multline*}
\frac{\partial}{\partial t^{1}}\frac{\partial L_{\alpha}}{\partial t^{2}%
}-\frac{\partial}{\partial t^{2}}\frac{\partial L_{\alpha}}{\partial t^{1}%
}=[[M_{\leq-1}^{(2)},M^{(1)}]_{\leq-1},L_{\alpha}]+[M_{\leq-1}^{(2)}%
,[M_{\leq-1}^{(1)},L_{\alpha}]]+\\
-[[M_{\leq-1}^{(1)},M^{(2)}]_{\leq-1},L_{\alpha}]-[M_{\leq-1}^{(1)}%
,[M_{\leq-1}^{(2)},L_{\alpha}]]
\end{multline*}

We have $[M_{\leq-1}^{(2)},[M_{\leq-1}^{(1)},L_{\alpha}]]-[M_{\leq-1}%
^{(1)},[M_{\leq-1}^{(2)},L_{\alpha}]]=[[M_{\leq-1}^{(2)},M_{\leq-1}%
^{(1)}],L_{\alpha}]$. Also $[M^{(2)},M^{(1)}]=0$ and therefore
\[
\lbrack M_{\leq-1}^{(2)},M^{(1)}]_{\leq-1}=-[M_{>-1}^{(2)},M^{(1)}]_{\leq
-1}=-[M_{>-1}^{(2)},M_{\leq-1}^{(1)}]_{\leq-1}=[M_{\leq-1}^{(1)},M_{>-1}%
^{(2)}]_{\leq-1}%
\]
Therefore $[M_{\leq-1}^{(2)},M^{(1)}]_{\leq-1}-[M_{\leq-1}^{(1)}%
,M^{(2)}]_{\leq-1}+[M_{\leq-1}^{(2)},M_{\leq-1}^{(1)}]=0$.
\end{proof}

\subsection{Bi-hamiltonian structure.}

\bigskip For any closed curve $x^{\alpha}=x^{\alpha}(s)$, $s\in S^{1}$ in
$\mathcal{U}$ restriction on $x(s)$ of flat connection corresponding to
$\{L_{\alpha}\}$ has form $L(s)=\partial/\partial s-\hbar^{-1}C$, where $C$ is
a function with values in $Mat(n,\mathbb{C)}$. We consider the following
brackets on the space of matrix-valued functions:
\begin{align*}
\{C_{i}^{j}(s),C_{k}^{l}(\widetilde{s})\}_{0} &  =\delta(s-\widetilde
{s})(\delta_{i}^{l}C_{k}^{j}(s)-\delta_{k}^{j}C_{i}^{l}(s))\\
\{C_{i}^{j}(s),C_{k}^{l}(\widetilde{s})\}_{1} &  =\delta^{\prime}%
(s-\widetilde{s})\delta_{k}^{j}\delta_{i}^{l}%
\end{align*}
Adaptation of standard arguments (see proof of proposition 1.8 from \cite{DS})
shows that $\{\cdot,\cdot\}_{0}-\hbar\{\cdot,\cdot\}_{1}$ is a Poisson
brackets for any $\hbar\in\mathbb{C}$. If $\{L_{\alpha}(t)\}$ satisfies
equations \ref{Lalpha(t)} then its restriction  on $x(s)$ satisfies the same
equation: $\partial L(s,t)/\partial t=[\varphi(b)_{\leq-1},L(s,t)]$. Recall
that an equation is called Hamiltonian if for any functional $f(C)$ we have
$df(C(t))/dt=\{f,H\}$. The same arguments as in propositions \ref{dress},
\ref{cordress} show that for a deformation $\widetilde{C}$ of $C$
corresponding to $L(t)|_{x(s)}$ there exists deformed dressing transformation
$T(s)\in Mat(n,\mathbb{C)[[\hbar]]}$ such that $T^{-1}(s)(\partial/\partial
s-\hbar^{-1}\widetilde{C})T(s)=\partial/\partial s+\sum_{k=-1}^{+\infty}%
\hbar^{k}h_{k}(s)$, $h_{k}(s)\in Diag$.

\begin{proposition}
\bigskip Equations \ref{Lalpha(t)} are Hamiltonian with respect to both
$\{\cdot,\cdot\}_{0}$ and $\{\cdot,\cdot\}_{1}$, the corresponding
Hamiltonians are $H_{b}=Tr\sum_{k=0}^{m-1}b_{k+1}\int_{x(s)}h_{k}ds$ and
$\widetilde{H}_{b}=Tr\sum_{k=-1}^{m-2}b_{k+2}\int_{x(s)}h_{k}ds$ respectively.
\end{proposition}

\begin{proof}
For a pair of functionals $f(C)$, $g(C)$:
\[
\{f,g\}:=\int\int\sum_{i,j,k,l}\frac{\delta f(C)}{\delta C_{i}^{j}(s)}%
\frac{\delta g(C)}{\delta C_{l}^{k}(\widetilde{s})}\{C_{i}^{j}(s),C_{k}%
^{l}(\widetilde{s})\}dsd\widetilde{s}%
\]
which gives
\begin{align*}
\{f,g\}_{0} &  =\int\sum_{i,j,k}\frac{\delta f(C)}{\delta C_{i}^{j}(s)}%
(\frac{\delta g(C)}{\delta C_{k}^{i}(s)}C_{k}^{j}(s)-\frac{\delta g(C)}{\delta
C_{j}^{k}(s)}C_{i}^{k}(s))ds\\
\{f,g\}_{1} &  =-\int\sum_{i,j}\frac{\delta f(C)}{\delta C_{i}^{j}(s)}%
\frac{\partial}{\partial s}(\frac{\delta g(C)}{\delta C_{j}^{i}(s)})ds
\end{align*}
Equation for $C(t)$ reads as $\frac{\partial C}{\partial t}=-[\varphi
(b)_{0},C]$ and we have also $[\varphi(b)_{0},C]=[\varphi(b)_{(-1)}%
,\frac{\partial}{\partial s}]$. We have $\frac{df(C(t))}{dt}=\int\sum
_{i,j}\frac{\delta f(C)}{\delta C_{i}^{j}(s)}\frac{\partial C_{j}^{i}%
}{\partial t}ds$. Therefore we must show that $\delta H_{b}/\delta C_{j}%
^{i}=-(\varphi(b)_{0})_{i}^{j}$ and that $\delta\widetilde{H}_{b}/\delta
C_{j}^{i}=-(\varphi(b)_{(-1)})_{i}^{j}$. Recall that $T^{-1}(\partial/\partial
s-\hbar^{-1}C)T=\partial/\partial s+h$, $h=\sum_{k=-1}^{+\infty}\hbar^{k}%
h_{k}$, $h_{k}\in Diag$ and that $b$ denotes $\sum_{k=1}^{m}\hbar^{-k}b_{k}%
$.We have
\begin{align*}
\frac{\delta H_{b}(C)}{\delta C_{j}^{i}(s)} &  =Tr(\sum_{k=0}^{m-1}%
b_{k+1}\frac{\delta h_{k}}{\delta C_{j}^{i}})=Tr(b\frac{\delta h}{\delta
C_{j}^{i}})_{(-1)}=Tr(b\frac{\delta(T^{-1}(\partial/\partial s-\hbar^{-1}%
C)T)}{\delta C_{j}^{i}})_{(-1)}=\\
&  =-Tr(bT^{-1}\frac{\hbar^{-1}\delta C}{\delta C_{j}^{i}}T)_{(-1)}%
(s)-Tr(b[T^{-1}\frac{\delta T}{\delta C_{j}^{i}},\partial/\partial
s+h])_{(-1)}(s)
\end{align*}
Notice that $[b,\partial/\partial s+h]=0$ since $b,h\in Diag$ and therefore
$Tr(b[T^{-1}\frac{\delta T}{\delta C_{j}^{i}},\partial/\partial
s+h])=Tr([b,T^{-1}\frac{\delta T}{\delta C_{j}^{i}}]\partial/\partial s+h)=0$.
We have also
\[
Tr(bT^{-1}\frac{\hbar^{-1}\delta C}{\delta C_{j}^{i}}T)_{(-1)}=Tr(TbT^{-1}%
\frac{\delta C}{\delta C_{j}^{i}})_{0}=(\varphi(b)_{0})_{i}^{j}%
\]
Therefore $\delta H_{b}/\delta C_{j}^{i}=-(\varphi(b)_{0})_{i}^{j}$. The proof
in the case of $\widetilde{H}_{b}$ is the same.
\end{proof}

\section{\bigskip Concluding remarks.}

It is important to understand how to generalize above hierarchies to the case
of Frobenius manifolds which are not semi-simple. It is necessary for example
for applications to theory of Gromov-Witten invariants where except for rare
cases like projective spaces and some other homogenous spaces Frobenius
manifolds defined by quantum cohomologies of Kahler manifolds are not semi-simple.

I planned initially to include a section describing the extension of
variations of $\frac{\infty}{2}-$ Hodge structures over higher times in the
framework of \cite{B1}. However this would increase significantly the volume
of this note. I plan to return to it in one of subsequent publications. I plan
also to write down some applications including formulas relating partition
functions of massive 2D topological field theories paired with gravity (see
\cite{W} for conjectures about properties of such partition functions) with
$\tau-$functions of above hierarchies.

\end{document}